\documentclass{article}
\usepackage[utf8]{inputenc}
\usepackage[abs]{overpic}
\usepackage{epsfig}
\usepackage{fancybox, ascmac}
\usepackage{multicol}
\usepackage{color}
\usepackage{mathrsfs}
\usepackage{amssymb}
\usepackage{bm}
\usepackage{amsmath}
\newtheorem{theorem}{Theorem}[section]

\newtheorem{example}[theorem]{Example}
\newtheorem{remark}[theorem]{Remark}

\title{Implementation of the Habegger--Lin decision algorithm}
\author{Yuka Kotorii
\thanks{Mathematics Program, Graduate School of Advanced Science and Engineering, Hiroshima University, 1-7-1 Kagamiyama Higashi-hiroshima City, Hiroshima 739-8521 Japan. E-mail: kotorii@hiroshima-u.ac.jp} 
\thanks{International Institute for Sustainability with Knotted Chiral Meta Matter (WPI-SKCM$^2$), Hiroshima University, 1-3-2 Kagamiyama Higashi-hiroshima City, Hiroshima 739-8521 Japan}
\thanks{RIKEN, interdisciplinary Theoretical and Mathematical Sciences Program, 2-1, Hirosawa, Wako, Saitama 351-0198, Japan}
\and 
Atsuhiko Mizusawa\thanks{Faculty of Science and Engineering, Waseda University, 3-4-1 Okubo, Shinjuku-ku, Tokyo 169-8555, Japan. E-mail: a\_mizusawa@aoni.waseda.jp}}

\begin{document}

\maketitle

\begin{abstract}
Habegger and Lin gave a classification of link-homotopy classes of links in terms of that of string links modulo certain group actions. 
As an application, they constructed an algorithm for determining whether given two links are link-homotopic.
In \cite{KM4}, we explicitly computed these group actions for the 4- and 5-component cases. Consequently, the Habegger--Lin algorithm can be effectively applied in these cases.
In this paper, we present an implementation of this algorithm, which is available at \cite{KMcode}, and exhibit new pairs of links that are not link-homotopic yet cannot be distinguished by Milnor's link-homotopy invariants, called $\overline{\mu}$-invariants.
\end{abstract}


2020 MATHEMATICS SUBJECT CLASSIFICATION: 57K10, 57M05.

\section{Introduction}
A {\it link} is an embedding of a disjoint union of circles into the 3-sphere $\mathbb{S}^3$ (see Figure~\ref{n-comp}, right). 
In this paper, we assume that a link is oriented and that its components are ordered.
Milnor \cite{Mil} introduced the notion of \textit{link-homotopy} on links, which is a restricted homotopy and a weaker equivalence relation than ambient isotopy, which we take as the usual equivalence relation.
The link-homotopy is the equivalence relation generated by ambient isotopies and self-crossing changes, where a self-crossing change is a crossing change between strands of the same component.
Milnor \cite{Mil} (see also \cite{Mil2}) gave a complete classification of link-homotopy classes of links with 3 or fewer components by numerical invariants called $\overline{\mu}$-invariants.
On the other hand, Levine \cite{Le} showed that the $\overline{\mu}$-invariants are not enough to classify 4-component links completely, and he \cite{Le2} refined indeterminacy of the $\overline{\mu}$-invariants and then classified 4-component links.

Subsequently, Habegger and Lin \cite{HL} provided a classification of link-homotopy classes for links with an arbitrary number of components.
To this end, they introduced a notion of a {\it string link} (see Figure~\ref{n-comp}, left), defined as an embedding of several intervals into $D^2\times [0,1]$ with fixed endpoints on the $D^2\times \{0\}$ and $D^2\times \{1\}$. 
They then showed that the set of link-homotopy classes of string links has a group structure and gave a Markov-type theorem for link-homotopy classes of string links. Namely, the closures of two string links are link-homotopic if and only if they are related by a finite sequence of string links such that any two consecutive elements differ by either a conjugation or a partial conjugation.
(Subsequently, Hughes \cite{Hu} showed that the partial conjugations generate the conjugations.) 
Moreover, based on the actions of the partial conjugations and the conjugations, they constructed an algorithm that determines whether the closures of two string links are link-homotopic.
Combined with the Markov-type theorem, this provided a method for deciding whether given two links are link-homotopic.
However, in general, these actions are not given explicitly, and thus the algorithm cannot be carried out in practice.

In \cite{KM4} and \cite{KM2}, by applying the clasper theory, independently introduced by Habiro \cite{Ha} and Gusarov \cite{G}, to 4- and 5-component string links, we explicitly calculated the actions of the partial conjugations and obtained classifications of link-homotopy classes of 4- and 5-component links. 
Graff \cite{Gra} also showed the classification for that of 4-component case independently by using similar techniques and for algebraically split 5-component links. 

In this paper, we present an implementation of this algorithm in Mathematica, which is available at \cite{KMcode}, and exhibit new pairs of links that are not link-homotopic, yet cannot be distinguished by $\overline{\mu}$-invariants.

\section{Preliminaries}
\subsection{Link presentation}

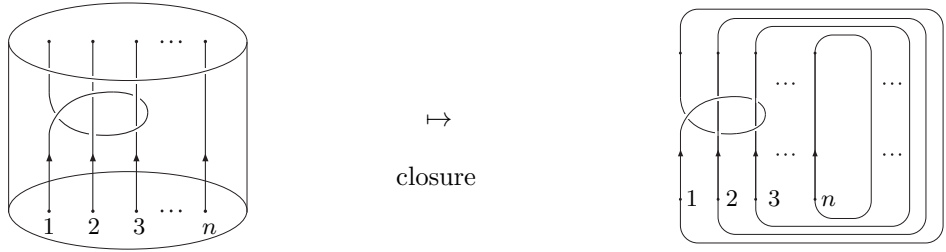
\begin{figure}[ht]
\begin{tabular}{ccc}
\begin{minipage}{0.4\hsize}
\centering
%\vspace{-3.4cm}
\includegraphics[width=90 pt]{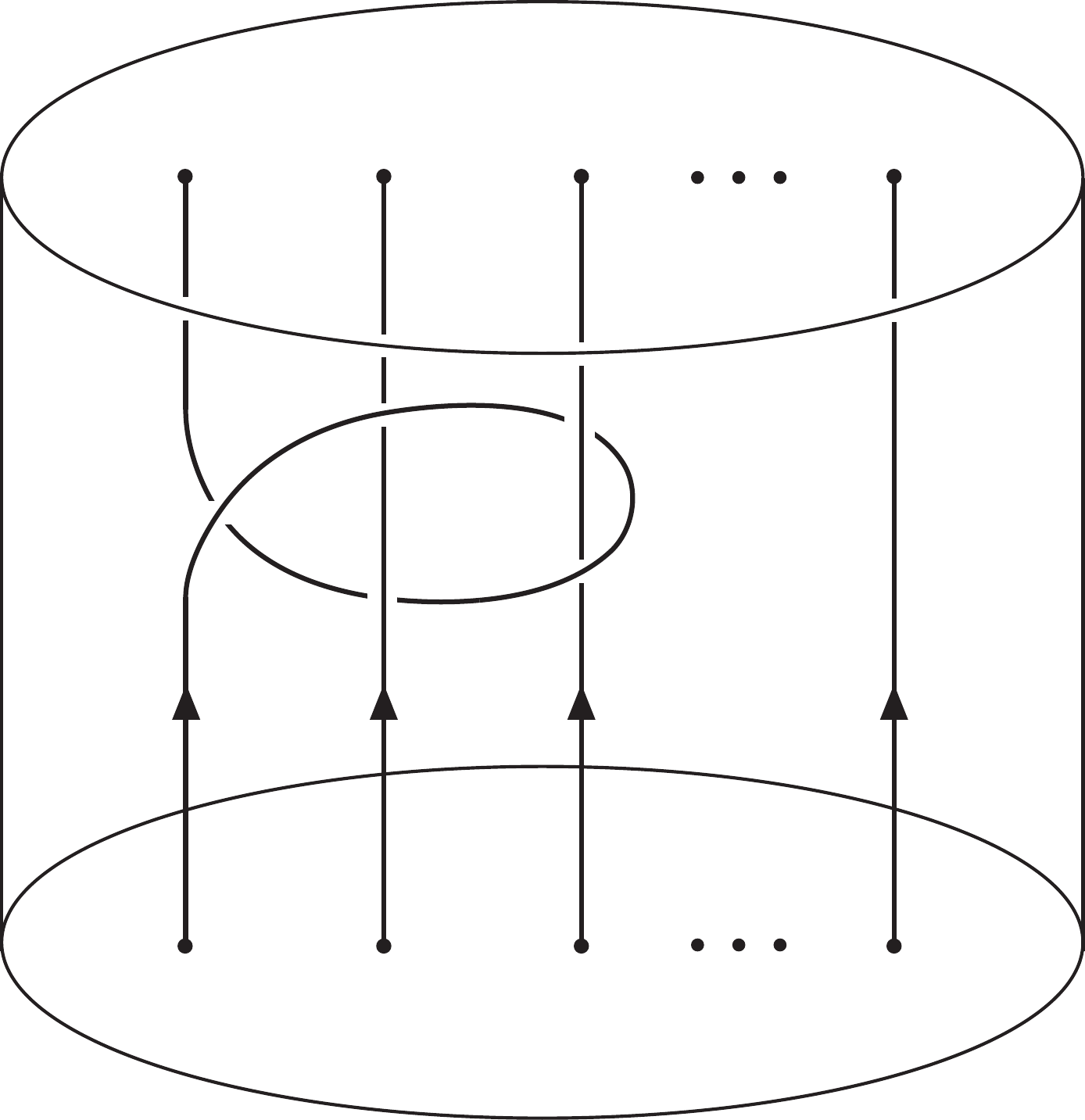}
\put(-77,6){{\small $1$}} \put(-60,6){\small $2$} \put(-43,6){\small $3$} 
\put(-17,6){\small $n$} 
\end{minipage}
\begin{minipage}{0.1\hsize}
$$\mapsto$$ 
\ \ closure
\end{minipage}
\begin{minipage}{0.5\hsize}
\centering
%\vspace{-3cm}
\includegraphics[width=100 pt]{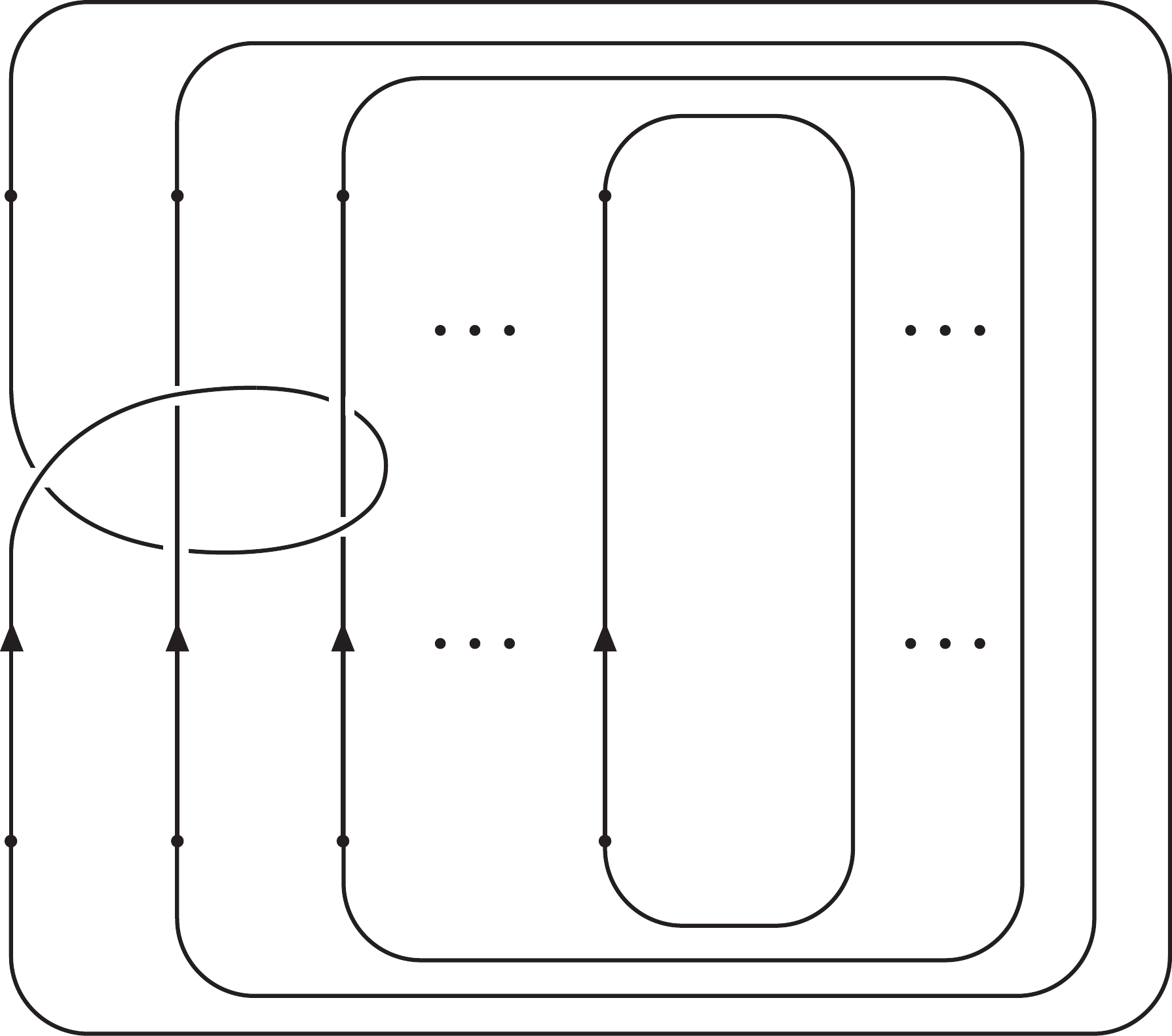}
\put(-97,14){{\small $1$}} \put(-82,14){\small $2$} \put(-66,14){\small $3$} 
\put(-46, 14){\small $n$} 
\end{minipage}
\end{tabular}
\caption{The closure of an $n$-component string link} \label{n-comp}
\end{figure}

Given an $n$-component string link, one obtains an $n$-component link as illustrated in Figure~\ref{n-comp}; this operation is called the {\it closure}.\\

The closure map is surjective up to ambient isotopy. 
Hence, any link can be represented as the closure of a string link. 
Indeed, for any $n$-component link $L$, by cutting each component at a point, we obtain an $n$-component string link whose closure is $L$. 
However, this representation is not unique.
For the link-homotopy case, Habegger and Lin \cite{HL} showed the Markov-type theorem as stated in Section~1.

\begin{figure}[ht]
\centering
\includegraphics[width=28pt, angle=90]{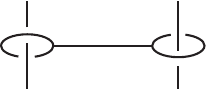}
\put(1,12){\large$:=$}
\put(75,12){\large$T_{1}$}
\hspace{0.7cm}
\includegraphics[width=28pt, angle=90]{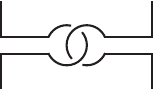}
%  \hspace{0.5cm} \mbox{\large$:=$}

\includegraphics[width=50pt]{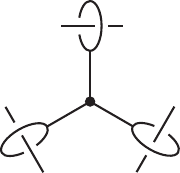}
\put(1,17){\large$:=$}
\put(75,12){\large$T_{2}$}
\hspace{0.7cm}
\includegraphics[width=50pt]{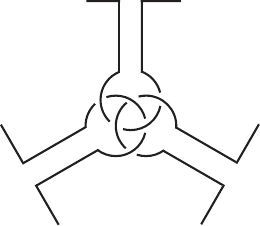} 
\vspace{.5cm} 

\includegraphics[bb=0 0 217 73, width=100pt]{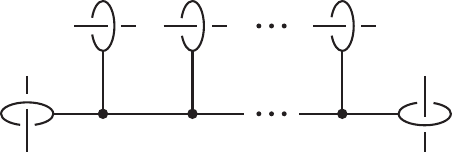}
\put(1,14){\large$:=$}
\put(105,42){\large$T_{n}$}
% \mbox{\large$:=$}
\hspace{0.6cm}
\includegraphics[width=100 pt]{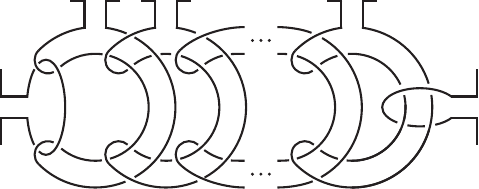}
\caption{Local structures $T_1, T_2$ and $T_n$} \label{clasper}
\end{figure}

\begin{figure}[ht]
\centering
\includegraphics[width=120 pt]{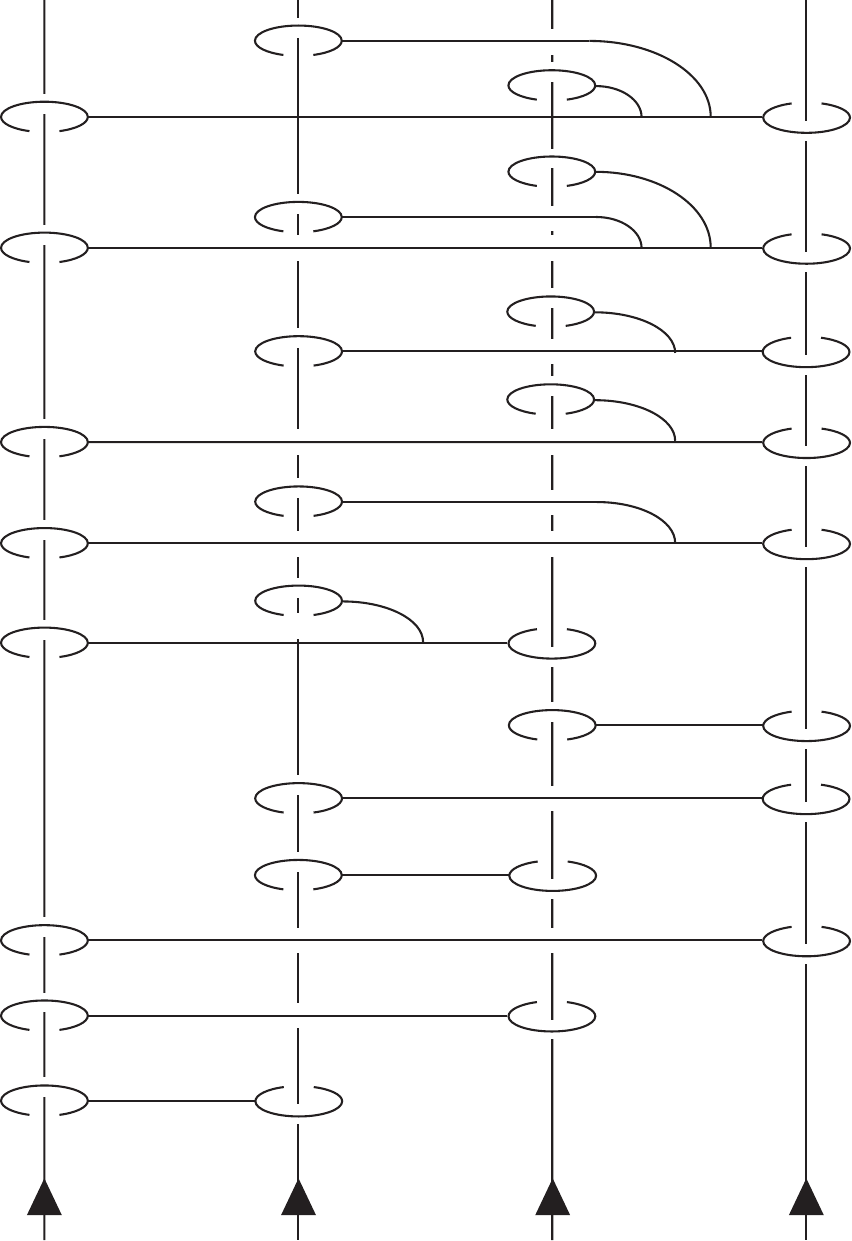}
\put(-116,-13){$1$}\put(-81,-13){$2$}
\put(-45,-13){$3$}\put(-9,-13){$4$}
\put(-69,18){$y_{12}$}\put(-33,29){$y_{13}$}
\put(3,40){$y_{14}$}\put(-33,51){$y_{23}$}
\put(3,60){$y_{24}$}\put(3,73){$y_{34}$}
\put(-33,83){$y_{123}$}\put(3,97){$y_{124}$}
\put(2,111){$y_{134}$}\put(2,124){$y_{124}$}
\put(2,139){$y_{1234}$}\put(2,157){$y_{1324}$}
\caption{The canonical form for 4-component string links} \label{4-comp}
\end{figure}

Meilhan and Yasuhara \cite{MY2} gave a canonical form for link-homotopy classes of string links with an arbitrary number of components. For example, in the case of 4-component string links, a link-homotopy class is determined by an integer sequence of length 12, 
\begin{align*}
& Y_1=(y_{12}, y_{13}, y_{14}, y_{23}, y_{24}, y_{34}), \\
& Y_2=(y_{123}, y_{124}, y_{134}, y_{234}), \\
& Y_3=(y_{1234}, y_{1324}), 
\end{align*}
where each entry $y_I$ (respectively, $-y_I$) represents the number of repetitions of the local structure $T_{|I|}$ as illustrated in Figure~\ref{clasper} (respectively, $T_{|I|}^{-1}$ which is obtained from $T_{|I|}$ by adding a positive half-twist \includegraphics[width=16pt]{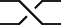} at the right-most strands) if $y_I$ is positive (respectively, negative) attached to the trivial string link, as depicted in Figure~\ref{4-comp}, according to the index sequence $I$, where $|I|$ denotes the length of the sequence $I$.\\

\begin{figure}[ht]
\centering
\includegraphics[width=270 pt]{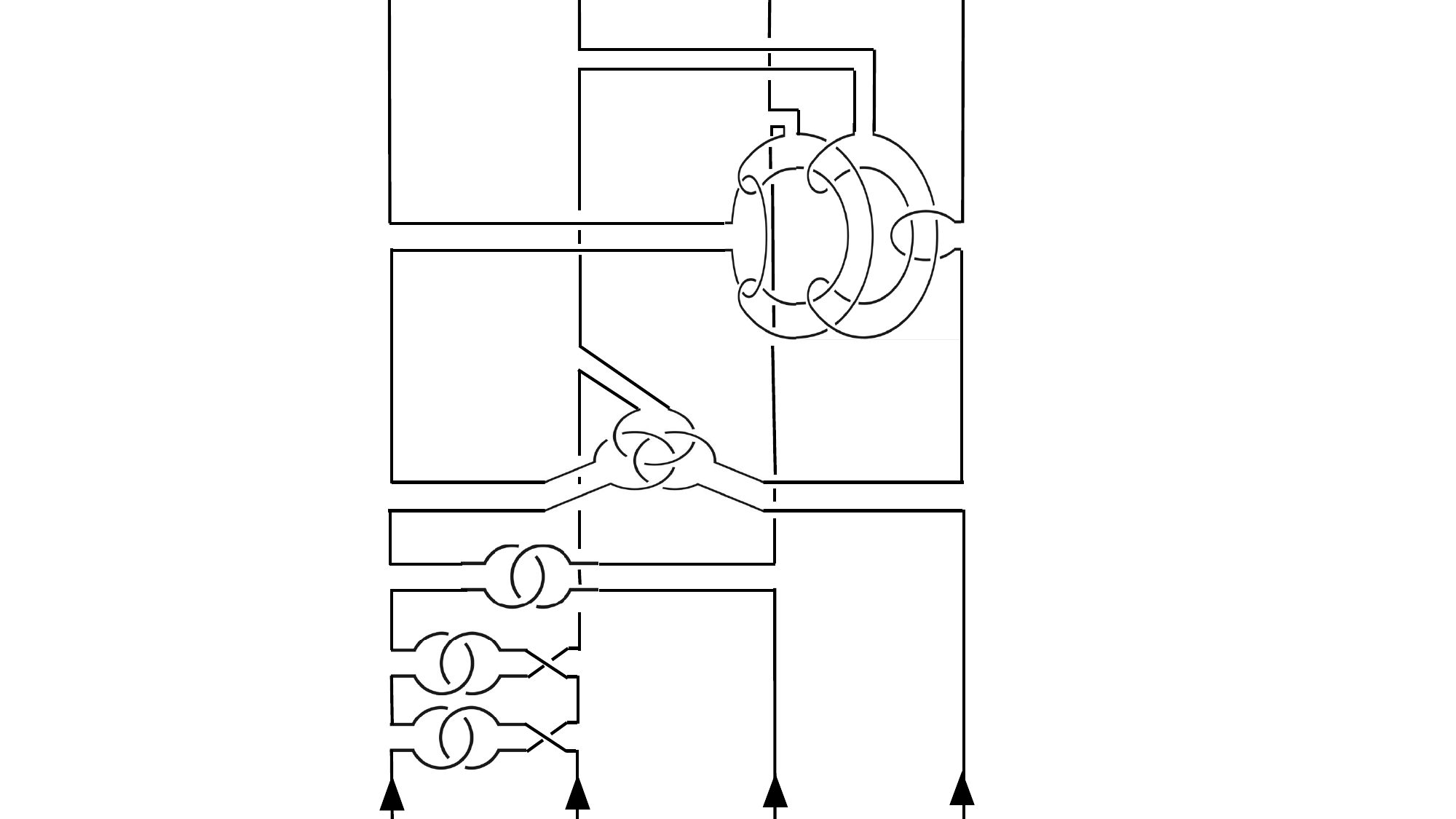}
\put(-200,-13){$1$}\put(-165,-13){$2$}
\put(-129,-13){$3$}\put(-93,-13){$4$}
\caption{A 4-component string link in the canonical form} \label{4-complink}
\end{figure}

For example, the integer sequence of length 12 
$$(Y_1, Y_2,Y_3)=(-2, 1, 0, 0, 0, 0, 0, 1, 0, 0, 0, 1),$$
represents the string link shown in Figure~\ref{4-complink}. 

Similarly, the link-homotopy class of a 5-component string link is determined by an integer sequence of length 36,
\begin{align*}
&Y_1=(y_{12}, y_{13}, y_{14}, y_{15}, y_{23}, y_{24}, y_{25}, y_{34}, y_{35}, y_{45}),\\ 
&Y_2=(y_{123}, y_{124}, y_{134}, y_{125}, y_{135}, y_{145}, y_{234}, y_{235}, y_{245}, y_{345}), \\ 
& Y_3=(y_{1234}, y_{1324}, y_{1235}, y_{1245}, y_{1325}, y_{1345}, y_{1425}, y_{1435}, y_{2345}, y_{2435}), \\ 
& Y_4=(y_{12345}, y_{12435}, y_{13245}, y_{13425}, y_{14235}, y_{14325}). 
\end{align*}
For the canonical form of 5-component string links, see Figure 11 in \cite{KM4}.

Using this canonical form, we obtained Theorems 3.2 and 3.6 in \cite{KM4}, which state that the set of link-homotopy classes of 4- (respectively, 5-) component links is in one-to-one correspondence with the set of ordered tuples of the above 12 (respectively, 36) integers modulo certain relations generated by the actions of partial conjugations (see \cite{KM4} for details).
This description is essentially equivalent to the classification by Habegger and Lin \cite{HL}.
Using the actions, they constructed the algorithm to determine whether given two links are link-homotopic.
However, they did not carry out explicit computations of the actions.
In the cases of 4- and 5-component links,
our description above enables us to execute the Habegger--Lin algorithm as described in the next subsection.

\subsection{Habegger--Lin algorithm}

Let $L$ and $L'$ be two $n$-component links represented, via the closure, by $(Y_1, Y_2, \dots, Y_{n-1})$ and $(Y'_1, Y'_2, \dots,$ $Y'_{n-1})$, respectively. 
Then the Habegger--Lin algorithm can be reformulated in the notation of \cite{KM4} as follows. 

We will change $Y_i$ to $Y'_i$ in order. Note that the partial conjugations do not change the numbers $y_{ij}$. 

\begin{enumerate}
\renewcommand{\labelenumi}{\arabic{enumi}.}
\item Check whether $Y_1 = Y'_1$. If so, proceed to Step~2; otherwise, $L$ and $L'$ are not link-homotopic.
\item Let $\mathcal{S}_1$ denote the set of partial conjugations. Find $\Psi_{1} \in \mathcal{S}_1$ such that $\Psi_{1}(Y_2) = Y'_2$ for $L$. If such a $\Psi_{1}$ exists, then replace
\[
L = (Y_1, Y_2, Y_3, Y_4, \dots, Y_{n-1})
\]
by
\begin{align*}
L_1 &= \Psi_1(L) \\
    &= (Y_1, \Psi_1(Y_2), \Psi_1(Y_3), \Psi_1(Y_4), \dots, \Psi_1(Y_{n-1})) \\
    &= (Y'_1, Y'_2, \Psi_1(Y_3), \Psi_1(Y_4), \dots, \Psi_1(Y_{n-1})),
\end{align*}
and proceed to Step~3; otherwise, $L$ and $L'$ are not link-homotopic.
\item Let $\mathcal{S}_2$ denote the set of partial conjugations that do not change $Y'_2$ of $L_1$. Find $\Psi_{2} \in \mathcal{S}_2$ such that $\Psi_2 \circ \Psi_1(Y_3) = Y'_3$ for $L_1$. If such a $\Psi_{2}$ exists, then replace $L_1$ by
\begin{align*}
L_2 &= \Psi_2(L_1) \\
    &= (Y'_1, Y'_2, \Psi_2 \circ \Psi_1(Y_3), \Psi_2 \circ \Psi_1(Y_4), \dots, \Psi_2 \circ \Psi_1(Y_{n-1})) \\
    &= (Y'_1, Y'_2, Y'_3, \Psi_2 \circ \Psi_1(Y_4), \dots, \Psi_2 \circ \Psi_1(Y_{n-1})),
\end{align*}
and proceed to Step~4; otherwise, $L$ and $L'$ are not link-homotopic.
\item Let $\mathcal{S}_3$ denote the set of partial conjugations that do not change $Y'_2$ and $Y'_3$ of $L_2$. Find $\Psi_{3} \in \mathcal{S}_3$ such that $\Psi_3 \circ \Psi_2 \circ \Psi_1(Y_4) = Y'_4$ for $L_2$. If such a $\Psi_{3}$ exists, then replace $L_2$ by
\begin{align*}
L_3 &= \Psi_3(L_2) \\
    &= (Y'_1, Y'_2, Y'_3, \Psi_3 \circ \Psi_2 \circ \Psi_1(Y_4), \dots, \Psi_3 \circ \Psi_2 \circ \Psi_1(Y_{n-1})) \\
    &= (Y'_1, Y'_2, Y'_3, Y'_4, \dots, \Psi_3 \circ \Psi_2 \circ \Psi_1(Y_{n-1})),
\end{align*}
and proceed to the next step; otherwise, $L$ and $L'$ are not link-homotopic.
\item Repeat this process until $Y_{n-1}$ is transformed into $Y'_{n-1}$. In that case,
\[
L' = \Psi_{n-1} \circ \cdots \circ \Psi_1(L),
\]
and hence $L$ and $L'$ are link-homotopic. Otherwise, if the procedure terminates before this stage, $L$ and $L'$ are not link-homotopic.
\end{enumerate}
Note that each step of the algorithm can be verified by solving the corresponding system of linear Diophantine equations.

\section{Main Results}

Based on \cite{KM4}, we provide an implementation of the Habegger--Lin algorithm for 4- and 5-component (string) links in Mathematica, and the code is available at \cite{KMcode}.
The input to this code consists of $(Y_i)_{i=1,2,\dots,n-1}$ and $(Y'_i)_{i=1,2,\dots,n-1}$ ($n=4$ or $5$), which represent two links via the closure.
The output determines whether they are link-homotopic or not.
As additional information, the step at which the algorithm determines that they are not link-homotopic is also provided. 
Note that, if the links are link-homotopic, the procedure proceeds to the final step. \\
%The values of the independent Milnor invariants are also output.

Using our code, we give three examples of pairs of links that are not link-homotopic, but all of their $\overline\mu$-invariants coincide.

\begin{remark}
Milnor's link-homotopy invariants $\overline{\mu}(I)$ are a family of link invariants indexed by non-repeating multi-indices $I$, taking values in cyclic groups.
In particular, $\overline{\mu}(ij)$ (for $i \neq j$) coincides with the linking number of the $i$th and $j$th components.
In \cite{Le2}, it is shown that, in the 4-component case,
the $\overline\mu$-invariants with the following indices $I$ determine all the others; hence, it suffices to examine them: 
\begin{align*}
& 12,13,14,23,24,34,\\
& 123,124,134,234,\\
& 1234,1324,
%& \overline{\mu}(12),\overline{\mu}(13),\overline{\mu}(14),\overline{\mu}(23),\overline{\mu}(24),\overline{\mu}(34),\\
%& \overline{\mu}(123), \overline{\mu}(124),\overline{\mu}(134),\overline{\mu}(234), \\
%& \overline{\mu}(1234), \overline{\mu}(1324)
\end{align*}
and in \cite{KM5} we proved that, in the 5-component case, it suffices to examine the $\overline\mu$-invariants with the following indices $I$: 
\begin{align*}
&12, 13, 14, 15, 23, 24, 25, 34, 35, 45, \\
&123, 124, 125, 134, 135, 145, 234, 235, 245, 345, \\
&1234, 1324, 1235, 1245, 1325, 1345, 1425, 1435, 2345, 2435, \\
&12345, 12435, 13245, 13425, 14235, 14325, \\
&21345, 21435, 31245, 31425, 41235, 41325.
%& \overline{\mu}(12),\overline{\mu}(13),\overline{\mu}(14),\overline{\mu}(15),\overline{\mu}(23),\\
%& \overline{\mu}(24),\overline{\mu}(25),\overline{\mu}(34),\overline{\mu}(35),\overline{\mu}(45), \\
%& \overline{\mu}(123),\overline{\mu}(124),\overline{\mu}(125),\overline{\mu}(134),\overline{\mu}(135), \\
%& \overline{\mu}(145), \overline{\mu}(234),\overline{\mu}(235),\overline{\mu}(245),\overline{\mu}(345), \\
%& \overline{\mu}(1234),\ \overline{\mu}(1324),\ \overline{\mu}(1235),\ \overline{\mu}(1245),\ 
%\overline{\mu}(1325),\\
%& \overline{\mu}(1345),\ \overline{\mu}(1425),\ \overline{\mu}(1435),\ 
%\overline{\mu}(2345),\ \overline{\mu}(2435) \\ 
%& \overline{\mu}(12345),\ \overline{\mu}(12435),\ \overline{\mu}(13245),\ \overline{\mu}(13425),\\
%& \overline{\mu}(14235),\ \overline{\mu}(14325),\ \overline{\mu}(21345),\ \overline{\mu}(21435),\\
%& \overline{\mu}(31245),\ \overline{\mu}(31425),\ \overline{\mu}(41235),\ \overline{\mu}(41325).
\end{align*}
\end{remark}

\begin{example}[Example~4.2. in \cite{KM4}]
Let $L$ and $L'$ be 5-component links, as illustrated in Figure~\ref{5-complink_exa1}, which are determined by
\begin{align*}
&Y_1=Y'_1 =(1, 0, 0, 0, 0, 0, 0, 0, 0, 0),\\ 
&Y_2=Y'_2=(0, 0, 0, 0, 0, 0, 0, 0, 0, 0), \\ 
&Y_3=Y'_3=(1, 0, 0, 0, 0, 0, 0, 0, 0, 0), \\ 
&Y_4=(0, 0, 1, 0, 0, 0), Y'_4=(0, 0, 0, 0, 0, 0). 
\end{align*}

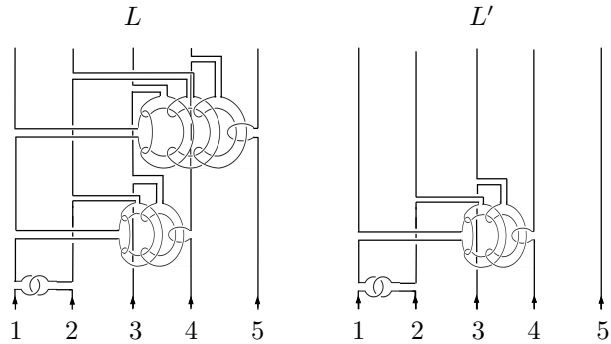
\begin{figure}[ht]
\centering
\includegraphics[width=230 pt]{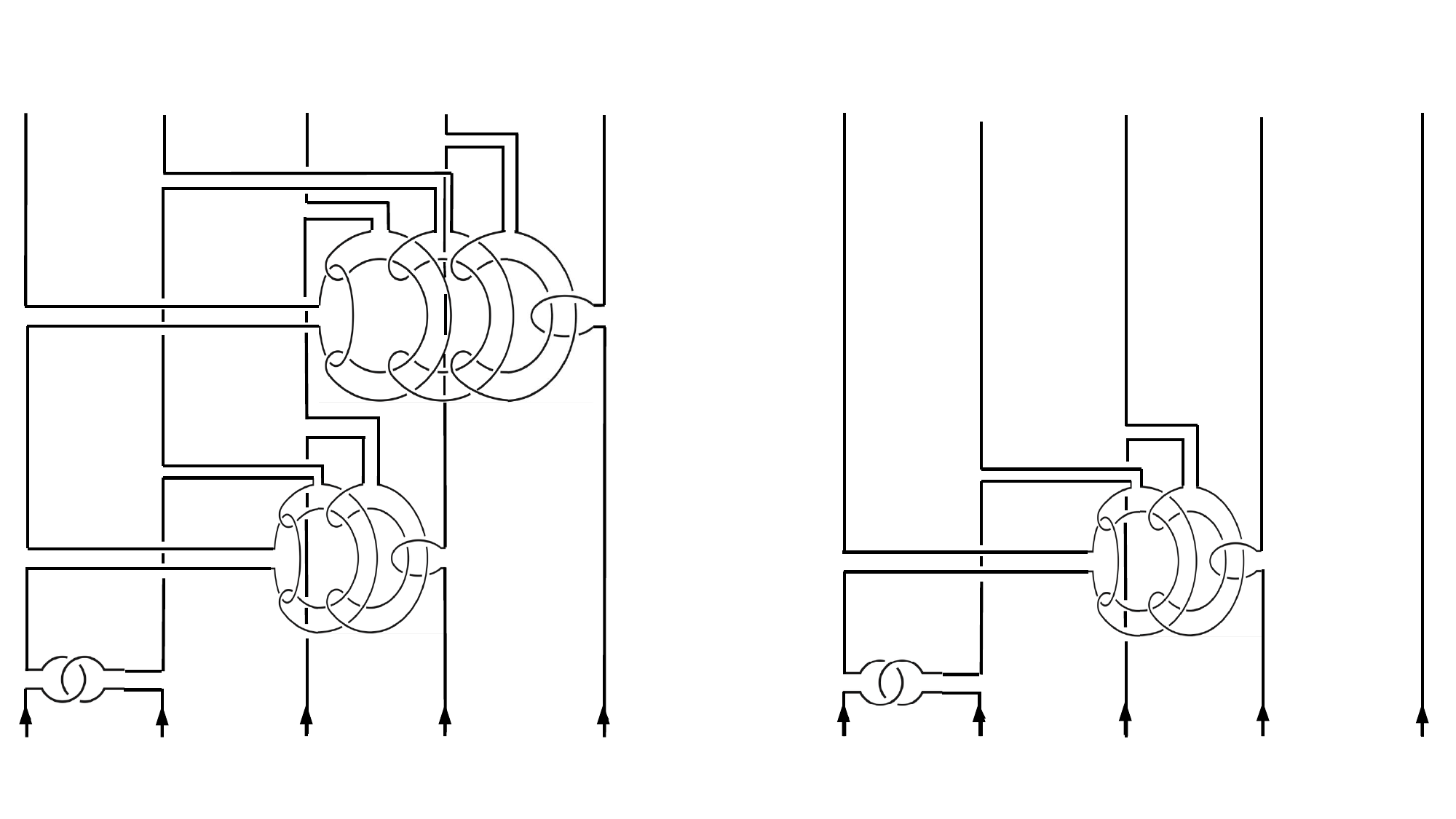}
\put(-185,120){$L$}\put(-55,120){$L'$}
\put(-228,2){$1$}\put(-207,2){$2$}
\put(-183,2){$3$}\put(-162,2){$4$}\put(-137,2){$5$}
\put(-99,2){$1$}\put(-77,2){$2$}
\put(-54,2){$3$}\put(-33,2){$4$}\put(-7,2){$5$}
\caption{5-component string links $L$ and $L'$} \label{5-complink_exa1}
\end{figure}

Using our implementation, it is shown (in Step~4) that $L$ and $L'$ are not link-homotopic, although they cannot be distinguished by Milnor's link-homotopy invariants. 
%Indeed, $\overline{\mu}(12)=1$ and the other Milnor invariants of length 2 are 0.
%Moreover, $\overline{\mu}(125)=\overline{\mu}(135)=\overline{\mu}(145)=\overline{\mu}(234)=\overline{\mu}(235)=\overline{\mu}(245)=\overline{\mu}(345)=\overline{\mu}(1345)=\overline{\mu}(1435)=\overline{\mu}(2345)=\overline{\mu}(2435)=0$ and the other Milnor invariants of length 3, 4 and 5 are $0 \mod 1$.
%$\overline{\mu}(123)=\overline{\mu}(124)=\overline{\mu}(125)=0 \mod 1$ and the other Milnor invariants of length 3 are $0$.
%Furthermore $\overline{\mu}(1345)=\overline{\mu}(1435)=\overline{\mu}(2345)=\overline{\mu}(2435)=0$ and the other Milnor invariants of length 4 are $0 \mod 1$.
Indeed, $\overline{\mu}(12)=1$ and the others are $0$ or $0 \mod 1$ for both links.
\end{example}

\begin{example}\label{exa2}
Let $L$ and $L'$ be 5-component links, as illustrated in Figure~\ref{5-complink_exa2}, which are determined by 
\begin{align*}
&Y_1=Y'_1 =(0, 0, 0, 0, 0, 0, 0, 0, 0, 0), \\ 
&Y_2=Y'_2=(1, 1, 0, 0, 0, 0, 0, 0, 0, 0), \\ 
&Y_3=(0, 0, 0, 0, 1, 0, 0, 0, 0, 0), Y'_3=(0, 0, 0, 0, 0, 0, 0, 0, 0, 0), \\ 
&Y_4=Y'_4=(0, 0, 0, 0, 0, 0).
\end{align*}
\begin{figure}[ht]
\centering
\includegraphics[width=230 pt]{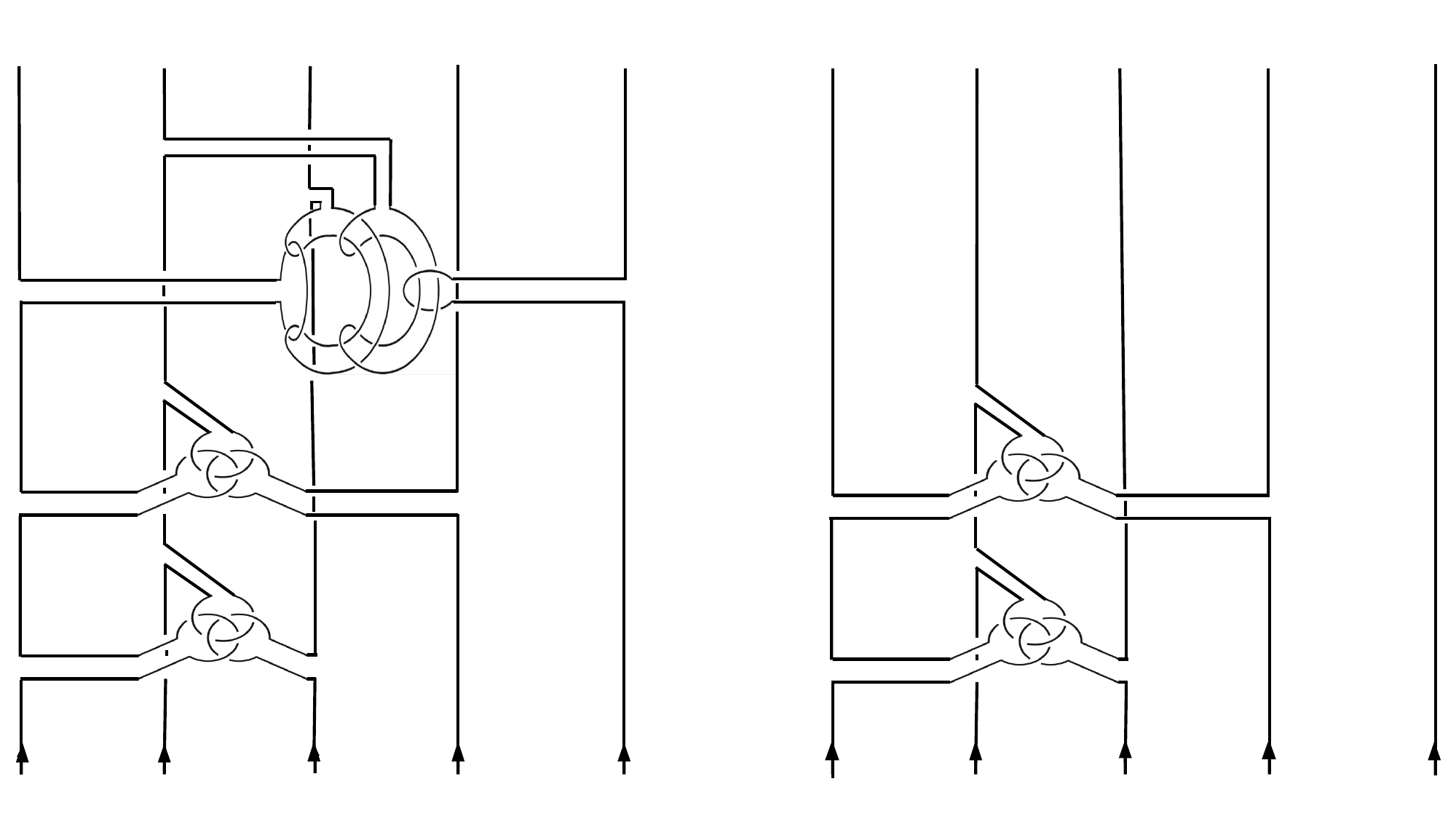}
\put(-185,125){$L$}\put(-55,125){$L'$}
\put(-229,-2){$1$}\put(-207,-2){$2$}
\put(-183,-2){$3$}\put(-161,-2){$4$}\put(-134,-2){$5$}
\put(-101,-2){$1$}\put(-78,-2){$2$}
\put(-54,-2){$3$}\put(-32,-2){$4$}\put(-5,-2){$5$}
\caption{5-component string links $L$ and $L'$} \label{5-complink_exa2}
\end{figure}
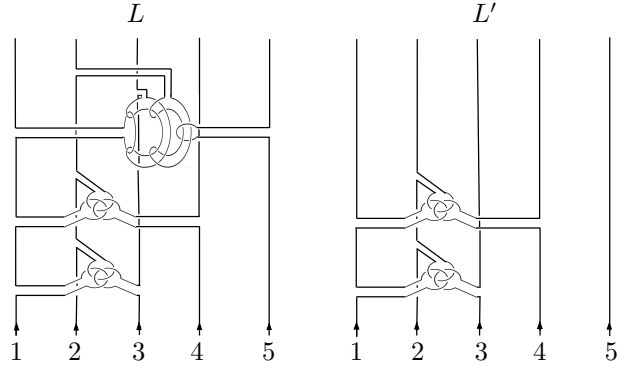

Using our implementation, it is shown (in Step~3) that $L$ and $L'$ are not link-homotopic, although they cannot be distinguished by Milnor's link-homotopy invariants. Indeed, we can check that all Milnor invariants of length 2 are $0$, $\overline{\mu}(123)=\overline{\mu}(124)=1$, and the other Milnor invariants of length 3 are $0$.
Moreover, $\overline{\mu}(1345)=\overline{\mu}(1435)=\overline{\mu}(2345)=\overline{\mu}(2435)=0$ and the other Milnor invariants of length 4 are $0 \mod 1$.
All Milnor invariants of length 5 are $0 \mod 1$ for $L$ and $L'$. 
\end{example}

\begin{example}\label{exa3}
Let $L$ and $L'$ be 5-component links determined by 
\begin{align*}
&Y_1=Y'_1 =(1, 1, 1, 1, 1, 1, 1,1, 1, 1), \\
&Y_2=(1, 1, 1, 1, 1, 1, 1, 1, 1, 1), Y'_2= (1, 1, 1, 1, 1, 1, 1, 1, 1, 0), \\ 
&Y_3=Y'_3=(1, 1, 1, 1, 1, 1, 1, 1, 1, 0), \\ 
&Y_4=Y'_4=(1, 1, 1, 1, 1, 1).
\end{align*}

Using our implementation, it is shown (in Step~2) that $L$ and $L'$ are not link-homotopic, although they cannot be distinguished by $\overline{\mu}$-invariants. 
Indeed, since all pairwise linking numbers are equal to $1$, all higher Milnor invariants are trivial modulo $1$ (see \cite{Mil2} for the definition of the Milnor invariants).
%Indeed, we can check that all $\overline{\mu}$-invariants of length 2 are $1$ and all $\overline{\mu}$-invariants of length 3, 4 and 5 are $0 \mod 1$ for both links. 
\end{example}

\begin{remark}
When $Y_1 = Y_1' = (1,1,\dots,1)$, the quantities $-y_{134}+y_{135}-y_{145}+y_{345}$ and $-y_{234}+y_{235}-y_{245}+y_{345}$
are invariant under Step 2 (see \cite{KM4} for details). 
In Example~\ref{exa3}, since these quantities take different values for $Y_2$ and $Y_2'$, it also follows that $L$ and $L'$ are not link-homotopic.
\end{remark}

\begin{remark}
In Example~\ref{exa2}, since the algorithm terminates at Step~3, for any $Y_4$ and $Y'_4$, the links $L$ and $L'$ are not link-homotopic. 
Similarly, in Example~\ref{exa3}, since the algorithm terminates at Step~2, for any $Y_3$, $Y_4$, $Y'_3$, and $Y'_4$, the links $L$ and $L'$ are not link-homotopic.
\end{remark}

\section{Acknowledgments}
The authors would like to thank Professor Takefumi Nosaka for his valuable comments and suggestions.
Y.K. is supported by JSPS KAKENHI, Grant-in-Aid for Scientific Research (C) Grant Number 25K07005, and by RIKEN iTHEMS Program and World Premier International Research Center Initiative  program, International Institute for Sustainability with Knotted Chiral Meta Matter (WPI-SKCM$^2$).

%\references

\end{document}